\documentclass[reqno,centertags]{amsart}
\usepackage{amsmath}
\usepackage{amssymb}
\usepackage{amsfonts}

\newtheorem{theorem}{Theorem}
\theoremstyle{plain}

\newtheorem{corollary}{Corollary}

\newtheorem{problem}{Problem}

\numberwithin{equation}{section}

\input{tcilatex}

\begin{document}
\title[]{Tubular surfaces whose Gauss map $N$ satisfies $\Delta ^{II}\boldsymbol{N}=\varLambda \boldsymbol{N}$}
\author{Hassan Al-Zoubi}
\address{Department of Mathematics, Al-Zaytoonah University of Jordan, P.O.
Box 130, Amman, Jordan 11733}
\email{dr.hassanz@zuj.edu.jo}
\date{}

\begin{abstract}
In this paper, we consider tubes in the Euclidean 3-space whose Gauss map $\boldsymbol{N}$ is of coordinate finite $II$-type, i.e., the  position vector $\boldsymbol{N}$ satisfies the relation $\Delta ^{II}\boldsymbol{N}=\varLambda \boldsymbol{N}$, where $\Delta ^{II}$ is the Laplace operator with respect to the second fundamental form $I$ of the surface and $\varLambda$ is a square matrix of order 3. We show that circular cylinders are the only class of surfaces mentioned above of coordinate finite $I$-type Gauss map.
\end{abstract}

\maketitle

\section{Introduction}
Euclidean immersions of finite Chen type regarding to the first fundamental form $I$ was introduced by B.-Y. Chen about four decades ago and it has been a topic of active research by many differential geometers since then. Many results in this area have been collected in \cite{C7}.

Let $\textbf{\textit{r}}$ be an isometric immersion of a surface $S$ in the 3-dimensional Euclidean space $\mathbb{E}^{3}$. We represent by $\Delta ^{I}$ the Laplacian operator of $S$ acting on the space of smooth functions $C^{\infty}(S)$. Then $S$ is called surface of finite $I$-type, if the position vector $\textbf{\textit{r}}$ of $S$ can be decomposed as a finite sum of eigenvectors of $\Delta ^{I}$ of $S$, that is
\begin{center}  \nonumber
$\textbf{\textit{r}}=\textbf{\textit{x}}_{0}+ \textbf{\textit{x}}_{1}+ \textbf{\textit{x}}_{2}+ ... + \textbf{\textit{x}}_{k},$
\end{center}%
where
\begin{center}  \nonumber
$\Delta ^{I}\textbf{\textit{x}}_{i}=\lambda _{i}\,\textbf{\textit{x}}_{i},\quad
i=1,... ,k,$
\end{center}
$\textbf{\textit{y}}_{0}$ is a constant vector and $\lambda _{1},\lambda_{2},... ,\lambda _{k}$ are eigenvalues of the operator $\Delta^{I}$.

Still now, the only known surfaces of finite type in the Euclidean 3-space are the minimal surfaces, the circular cylinders and the spheres. So in [18] B.-Y. Chen mentions the following problem
\begin{problem}
\label{(1)}Determine all surfaces of finite Chen $I$-type in $\mathbb{E}^{3}$.
\end{problem}

In order to provide an answer to the above problem, important families of surfaces were studied by different authors by proving that finite type ruled
surfaces, finite type quadrics, finite type tubes, finite type cyclides of Dupin and finite type spiral surfaces are surfaces of the only known
examples in $\mathbb{E}^{3}$. However, for another classical families of surfaces, such as surfaces of revolution, translation surfaces as well as helicoidal surfaces, the classification of its finite type surfaces is not known yet.

Following this ideas, in [17] B.-Y. Chen and P. Piccini studied submanifolds of which their Gauss map are of finite type with respect to the first fundamental form, so, specifically for the surfaces in $\mathbb{E}^{3}$, we are led to the following problem

\begin{problem}
\label{(2)}Classify all surfaces in $\mathbb{E}^{3}$ with finite type Gauss map corresponding to the first fundamental form $I$.
\end{problem}

Later in \cite{C3} O. Garay generalized Chen's condition studied surfaces for which surfaces in $\mathbb{E}^{3}$ satisfy the condition $\Delta^{I}\boldsymbol{r}= A\boldsymbol{r}$ $(\ddag)$ where $A \in\mathbb{Re}^{3\times3}$ or $\Delta^{I}\boldsymbol{N}= A\boldsymbol{N}$ where $\boldsymbol{N}$ is the Gauss map of $S$. It was shown that a surface $S$ in $\mathbb{E}^{3}$ satisfies $(\ddag)$ if and only if it is an open part of a minimal surface, a sphere, or a circular cylinder. Surfaces satisfying $(\ddag)$ are said to be of coordinate finite type, meanwhile Surfaces satisfying $\Delta^{I}\boldsymbol{N}= A\boldsymbol{N}$ are said to be of coordinate finite type Gauss map.

In the framework of surfaces of finite type Stamatakis and Al-Zoubi in \cite{S1} defined the notion of surfaces of finite type associated with the 2nd or 3rd fundamental form. For results concerning this see (\cite{A3}, \cite{A11}, \cite{A6}, \cite{A4}, \cite{A9}, \cite{A14}, \cite{A5}, \cite{A15}, \cite{A17}).

A general study which seems to be interesting can be made by investigating families of surfaces in $\mathbb{E}^{3}$ whose position vector $\textbf{\textit{r}}$ satisfies the relation
\begin{equation}
\Delta ^{J}\textbf{\textit{r}}=A\textbf{\textit{r}}, \ \ \ J = II, III,  \label{5}
\end{equation}%
where $A\in \mathbb{Re}^{3\times 3}$. Surfaces whose position vector $\textbf{\textit{r}}$ satisfies the above relation are said to be of coordinate finite $J$-type.

In \cite{A1} H. Al-Zoubi and S. Stamatakis, studied the quadric surfaces satisfying
\begin{equation}
\Delta ^{III}\textbf{\textit{r}}=A\textbf{\textit{r}}.   \label{6}
\end{equation}

The same authors in \cite{S2} classified the class of surfaces of revolution satisfying condition (\ref{6}). Later in \cite{A2}, the translation surfaces were studied, and it was proved that Sherk's surface is the only translation surface satisfying (\ref{6}).
B. Senoussi and M. Bekkar, in \cite{S0} classified the Helicoidal surfaces with $\Delta^{J}r = Ar, J = I, II, III$.

On the other hand, H. Al-Zoubi, T. Hamadneh in \cite{A14} studied the surfaces of revolution satisfying
\begin{center}  \nonumber
$\Delta ^{II}\textbf{\textit{r}}=A\textbf{\textit{r}}.$   \label{7}
\end{center}
where $A\in \mathbb{Re}^{3\times 3}$.

Also, the notion of surfaces of finite type can be applied to any smooth map not necessary the position vector of the surface, for example, its Gauss map. For more details, the reader can be referred to (\cite{A10}, \cite{A12}, \cite{A8}, \cite{A13}, \cite{B2}, \cite{B3}, \cite{B4}, \cite{B5}, \cite{B6}, \cite{S00}).

In this research, we will focus on surfaces of finite $II$-type. Firstly, we will define the second differential parameter of Beltrami regarded the 2nd fundamental form of $S$. Then, we complete our study by proving our main result.


\section{Preliminaries}
We consider a $C^{m}$-surface $S$, $m\geq3$, in the Euclidean 3-space $\mathbb{E}^{3}$ defined by an injective $C^{m} $-immersion $\boldsymbol{r} = \boldsymbol{r}(v^1,v^2)$ on a region $W: = U \times \mathbb{R} \,\,\,(U\subset \mathbb{R}$ open interval) of $\mathbb{R}^{2}$ with nonvanishing Gauss curvature. Denoting by
\begin{center}  \nonumber
$I = g_{st}\,d v^s d v^t, \quad II = b_{st}\,d v^s d v^t$,
\end{center}
\begin{center}  \nonumber
$III = e_{st}\,d v^s d v^t, \quad s,t = 1,2$,
\end{center}
the 1st, 2nd and 3rd fundamental forms of $S$ respectively. For any two sufficiently differentiable functions $p(v^{1}, v^{2})$ and $q(v^{1}, v^{2})$ on $S$, the 1st differential parameter of Beltrami according to the fundamental form $J$ is defined by \cite{A7}
\begin{equation}  \label{4.1}
\nabla^{J}(p,q):=a^{st}p_{/s}q_{/t},\ \ \ J = I, II, III
\end{equation}
where $p_{/s}: =\frac{\partial p}{\partial v^{s}}$, and $(a^{st})$ stand as the components of the inverse tensor of $(g_{st}), (b_{st})$ and $(e_{st})$ respectively.
The 2nd differential parameter of Beltrami of $S$ is defined by \cite{A7}
\begin{equation}  \label{4.2}
\triangle^{J}p:=-a^{st}\nabla^{J}_{s}p_{/t}, \ \ \ J = I, II, III,
\end{equation}
where $\nabla^{J}_{s}$ is the covariant derivative in the $v^{s}$ direction corresponding to the fundamental form $J$.

After a long computations, one can find that \cite{S1} 
\begin{equation}  \label{4.11}
\Delta^{II}\boldsymbol{N}=\frac{1}{2K}grad^{I}(K)+ 2H\boldsymbol{N}.
\end{equation}
where $K, H$ are the Gauss and mean curvature of $S$ respectively.
Next, we want to investigate the tubular surfaces whose Gauss map $\boldsymbol{N}$ satisfies the following relation
\begin{equation}   \label{4.22}
\Delta ^{II}\textbf{\textit{N}}=A\textbf{\textit{N}}.
\end{equation}

\section{Tubes in $\mathbb{E}^{3}$}
Let $\ell: \textbf{\textit{a}} =\textbf{\textit{a}}(\emph{u})$, $\mathit{u}\epsilon (a,b)$ be a regular unit speed curve of finite length. We denote by ${\textbf{\textit{t}}, \textbf{\textit{h}}, \textbf{\textit{b}}}$ the Frenet frame of the curve $\textbf{\textit{a}}$ and with $\kappa >0$ its curvature. Then a regular parametric representation of a tubular surface $\mathfrak{F}$ of radius $r$ satisfies $0 < r < min\frac{1}{|\kappa|}$ is given by

\begin{center}  \nonumber
$\mathfrak{F}:\textbf{\textit{x}}(u,\phi)= \textbf{\textit{a}}+ r \cos\phi\textbf{\textit{h}}+ r \sin\phi\textbf{\textit{b}}.$
\end{center}

One can find that fundamental forms $I$ and $II$ of $\mathfrak{F}$ are given by

\begin{center}  \nonumber
$I = \big(\delta^{2}+r^{2}\tau^{2}\big)du^{2} + 2r^{2}\tau dud\phi+r^{2}d\phi^{2},$
\end{center}
\begin{center}  \nonumber
$II = \big(-\kappa\delta\cos\phi+r\tau^{2}\big)du^{2} + 2r\tau dud\phi +rd\phi^{2},$
\end{center}
where $\delta: = (1-r\kappa\cos\phi)$ and $\tau$ is the torsion of the curve $\textbf{\textit{a}}$. The Gauss curvature of $\mathfrak{F}$ is given by

\begin{center}    \nonumber
$K =-\frac{\kappa\cos\phi}{r\delta}$.
\end{center}

In the following we consider an open portion of the tube surface where $K \neq 0 (\Leftrightarrow \cos\phi \neq 0)$. From (\ref{4.2}) the Laplace operator $\Delta ^{II}$ of $\mathfrak{F}$ is (see also \cite{A3}, formula (9))

\begin{eqnarray}  \label{eq8}
\Delta ^{II} &=&\frac{1}{\kappa\delta \cos \phi }\Bigg[\frac{\partial^{2}}{\partial u^{2}}-2\tau \frac{\partial ^{2}}{\partial u\partial \phi } \\  \nonumber
& &+\Big(\tau ^{2}-\frac{\kappa\delta \cos \phi}{r}\Big)\frac{\partial^{2}}{\partial \phi^{2}}+\frac{(1-2\delta)\beta }{2\kappa\delta \cos \phi }\frac{\partial }{\partial u}  \\  \nonumber
& &+\left(-\tau \acute{} +\frac{\tau \beta(2\delta-1) }{2\kappa\delta \cos \phi }+\frac{\kappa(2\delta-1)\sin\phi}{2r}\right) \frac{\partial }{\partial \phi }\Bigg] ,
 \nonumber
\end{eqnarray}
where $\beta: =\kappa \acute{}\cos \phi +\kappa \tau \sin \phi $ and $\acute{}:=\frac{d}{du}$.

The Gauss map $\textbf{\textit{N}}(u,\phi)$ of $\mathfrak{F}$ has parametric representation
\begin{equation}  \label{eq81}
\textbf{\textit{N}}(u,\phi)= -\cos \phi \textbf{\textit{h}} - \sin \phi \textbf{\textit{b}}.
\end{equation}

Using (\ref{4.11}) and (\ref{eq8}), we obtain
\begin{eqnarray}    \label{eq831} \nonumber
\Delta ^{II} \textbf{\textit{N}}&=&\frac{\beta}{2\kappa\delta^{2}\cos \phi}\textbf{\textit{t}}+ \bigg(\frac{\sin^{2}\phi}{2r\delta \cos\phi}+\frac{\cos \phi}{r\delta}-\frac{2\cos \phi}{r}\bigg)\textbf{\textit{h}} \\
& &+\frac{(1-4\delta)\sin\phi}{2r\delta}\textbf{\textit{b}}.
\end{eqnarray}

Let $(t_{i}, h_{i}, b_{i})$ and $N_{i}, i = 1,2,3$ be the coordinate functions of $\boldsymbol{t}, \boldsymbol{h}, \boldsymbol{b}, \boldsymbol{N}$ respectively, and let $a _{ij}$ for $i,j=1,2,3$ be the entries of the matrix $A$. By virtue of (\ref{4.22}) and (\ref{eq831}) one can find

\begin{eqnarray}    \label{eq834} \nonumber
&&\frac{\beta}{2\kappa\delta^{2}\cos \phi}t_{i}+ \bigg(\frac{\sin^{2}\phi}{2r\delta \cos\phi}+\frac{\cos \phi}{r\delta}-\frac{2\cos \phi}{r}\bigg)h_{i} \\  \nonumber
&&+\frac{(1-4\delta)\sin\phi}{2r\delta}b_{i}= \\  \nonumber
&&-a_{i1}(h_{1}\cos \phi + b_{1}\sin \phi) -a_{i2}(h_{2}\cos \phi + b_{2}\sin \phi) \\
&&-a_{i3}(h_{3}\cos \phi + b_{3}\sin \phi).
\end{eqnarray}

We rewrite the last equation as follows
\begin{eqnarray}    \label{eq832}
&&r\beta t_{i}+ \kappa\delta\big(1+(1-4\delta)\cos^{2} \phi\big)h_{i} \\  \nonumber
&&+\kappa\delta(1-4\delta)\cos \phi\sin\phi b_{i} \\ \nonumber
&&+2\kappa\delta^{2}\cos \phi\big(a_{i1}(h_{1}\cos \phi + b_{1}\sin \phi) \\  \nonumber
&&+a_{i2}(h_{2}\cos \phi + b_{2}\sin \phi) \\  \nonumber
&&+a_{i3}(h_{3}\cos \phi + b_{3}\sin \phi)\big)=0.
\end{eqnarray}

We have the following two cases

Case I. $\beta \equiv 0$. Then $\kappa' = 0$ and $\kappa\tau = 0$. Thus $\tau = 0$ and $\kappa$ = const. $\neq 0$, therefore the curve $\textbf{\textit{a}}$ is a plane circle and so, $\mathfrak{F}$ is an anchor ring.
In this case, the first fundamental form becomes



\begin{center}  \nonumber
$I=\delta^{2}du^{2}+r^{2}d\phi ^{2}$,
\end{center}%
while the second is

\begin{center}  \nonumber
$II= -\kappa\delta\cos\phi du^{2} +rd\phi^{2}.$
\end{center}

Hence, the Laplacian is given by  

\begin{equation}  \label{eq10}
\Delta^{II}= \frac{1}{\kappa \delta \cos \phi}\frac{\partial^{2}}{\partial u^{2}}-\frac{1}{r}\frac{\partial^{2}}{\partial\phi^{2}} +\frac{(2\delta-1)\sin \phi}{2r\delta \cos \phi}\frac{\partial}{\partial \phi}.
\end{equation}

For the Gauss map $\boldsymbol{N}$ of $\mathfrak{F}$, we have
\begin{equation*}
\boldsymbol{N}=\{-\cos u \cos \phi , -\sin u \cos\phi , -\sin \phi\}.
\end{equation*}

Let $\boldsymbol{N}=(N_{1}, N_{2}, N_{3})$. By virtue of (\ref{eq10}) we find

\begin{equation}  \label{eq11}
\Delta^{II}N_{1} = \Bigg[\frac{1}{\kappa\delta} -\frac{\cos\phi}{r}+\frac{(2\delta-1)\sin^{2}\phi}{2r\delta\cos\phi}\Bigg]\cos u,
\end{equation}

\begin{equation}  \label{eq12}
\Delta^{II}N_{2} =\Bigg[\frac{1}{\kappa\delta} -\frac{\cos\phi}{r}+\frac{(2\delta-1)\sin^{2}\phi}{2r\delta\cos\phi}\Bigg]\sin u,
\end{equation}

\begin{equation}  \label{eq13}
\Delta^{II}N_{3} = -\frac{\sin \phi}{2r\delta}(4\delta -1).
\end{equation}

Let $a _{ij}, i,j=1,2,3$ be the entries of the matrix $A$. From (\ref{4.22}) and on account of (\ref{eq11}), (\ref{eq12}) and (\ref{eq13}) we get

\begin{eqnarray}  \label{eq14}
& &\Bigg[\frac{1}{\kappa\delta} -\frac{\cos\phi}{r}+\frac{(2\delta-1)\sin^{2}\phi}{2r\delta\cos\phi}\Bigg]\cos u= \\ \nonumber
& &-a _{11}\cos u \cos \phi-a _{12}\sin u \cos\phi- a _{13}\sin \phi,
\end{eqnarray}

\begin{eqnarray}  \label{eq15}
& &\Bigg[\frac{1}{\kappa\delta} -\frac{\cos\phi}{r}+\frac{(2\delta-1)\sin^{2}\phi}{2r\delta\cos\phi}\Bigg]\sin u=\\ \nonumber
& &-a _{21}\cos u \cos \phi- a _{22}\sin u \cos\phi- a _{23}\sin \phi,
\end{eqnarray}%

\begin{eqnarray}  \label{eq16}
& &-\frac{\sin \phi}{2r\delta}(4\delta -1)=\\ \nonumber
& &- a _{31}\cos u \cos \phi- a _{32}\sin u \cos\phi- a _{33}\sin \phi.
\end{eqnarray}

From the last equation it is easily verified that $a_{31}=a_{32}=0$. Hence relation (\ref{eq16}) becomes
\begin{equation}  \label{eq17}
\frac{\sin \phi}{2r\delta}(4\delta -1)=  a _{33}\sin \phi,
\end{equation}
then it follows that
\begin{equation}  \label{eq18}
\frac{4\delta -1}{2r\delta}=  a _{33}.
\end{equation}

In this case, the above relation is valid only when $a_{33}$ is function of the parameters $u$ and $\phi$, Which contradicts our assumption. So we have
\begin{corollary}  \label{C1.2}
Every anchor ring in the Euclidean 3-space is of coordinate infinite type Gauss map regarding the second fundamental form.
\end{corollary}

Case II. $\beta \neq 0$. Recalling equations (\ref{eq832}) for $i = 1, 2, 3$, we obtain that these equations are polynomials in $\cos \phi,\sin\phi$ with coefficients functions of $u$. Therefore we conclude that
\begin{equation}  \label{eq19}
r\beta t_{i}=0,
\end{equation}
\begin{equation}  \label{eq20}
\kappa\delta\bigg(1+(1-4\delta)\cos^{2} \phi\bigg) h_{i}= 0,
\end{equation}
\begin{equation}  \label{eq21}
\kappa\delta(1-4\delta)\cos \phi\sin\phi b_{i}= 0,
\end{equation}
\begin{equation}  \label{eq22}
2\kappa\delta^{2}\cos^{2} \phi\big(a_{i1}h_{1}+a_{i2}h_{2}+a_{i3}h_{3}\big)=0,
\end{equation}
\begin{equation}  \label{eq23}
2\kappa\delta^{2}\cos \phi\sin\phi\big(a_{i1}b_{1}+a_{i2}b_{2}+a_{i3}b_{3}\big)=0.
\end{equation}

Since (\ref{eq19}) is valid for $i = 1, 2, 3$, so in vector notation it follows that
\begin{equation}  \label{eq24}
r\beta \boldsymbol{t}=\boldsymbol{0},
\end{equation}
from which we must have $\beta = 0$ and thus $\mathfrak{F}$ is an anchor ring. So, according to Corollary 1, we have
\begin{theorem}
All tubes in the Euclidean 3-space is of coordinate infinite type Gauss map regarding the second fundamental form.
\end{theorem}


\begin{thebibliography}{00}
\bibitem{A3} H. Al-Zoubi, Tubes of finite $II$-type in the Euclidean 3-space, WSEAS Trans. Math. \textbf{17} (2018), 1-5.
\bibitem{A10} H. Al-Zoubi, On the Gauss map of quadric surfaces, arXiv: 1905.00962v1, (2019).
\bibitem{A11} H. Al-Zoubi, W. Al Mashaleh, Surfaces of finite type with respect to the third fundamental form, IEEE Jordan International Joint Conference on Electrical Engineering and Information Technology (JEEIT), Amman, April 9-11, (2019).
\bibitem{A12} H. Al-Zoubi, M. Al-Sabbagh, Anchor rings of finite type Gauss map in the Euclidean 3-space, Int. J. Math. and Comput. Methods \textbf{5} (2020), 9-13.
\bibitem{A6} H. Al-Zoubi, M. Al-Sabbagh, S. Stamatakis, On surfaces of finite Chen $III$-type, Bull. Belgian Math. Soc. \textbf{26} (2019), 177-187.
\bibitem{A8} H. Al-Zoubi, H. Alzaareer, T. Hamadneh, M. Al Rawajbeh, Tubes of coordinate finite type Gauss map in the Euclidean 3-space, Indian J. Math. \textbf{62} (2020), 171-182.
\bibitem{A4} H. Al-Zoubi, S. Al-Zu’bi, S. Stamatakis and H. Almimi, Ruled surfaces of finite Chen-type. J. Geom. Graph. \textbf{22} (2018), 15-20.
\bibitem{A7} H. Al-Zoubi, A. Dababneh, M. Al-Sabbagh, Ruled surfaces of finite $II$-type, WSEAS Trans. Math. \textbf{18} (2019), 1-5.
\bibitem{A9} H. Al-Zoubi, T. Hamadneh, Surfaces of revolution of finite $III$-type, arXiv: 1907.12390v2, Oct (2019).
\bibitem{A13} H. Al-Zoubi, T. Hamadneh, Quadric surfaces of coordinate finite type Gauss map, arXiv 2006.04529v1, May (2020).
\bibitem{A14} H. Al-Zoubi, T. Hamadneh, Surfaces of coordinate finite $II$-type, arXiv: 2005.05120v1, May (2020).
\bibitem{A5} H. Al-Zoubi, K. M. Jaber, S. Stamatakis, Tubes of finite Chen-type, Comm. Korean Math. Soc. \textbf{33} (2018), 581-590.
\bibitem{A2} H. Al-Zoubi, S. Stamatakis, W. Al Mashaleh and M. Awadallah, Translation surfaces of coordinate finite type, Indian J. Math. \textbf{59} (2017), 227-241.
\bibitem{A1} H. Al-Zoubi, S. Stamatakis, Ruled and quadric surfaces satisfying $\triangle ^{III}\mathbf{x}=A\mathbf{x}$, J. Geom. Graph. \textbf{20} (2016), 147-157.
\bibitem{A15} H. Al-Zoubi, T. Hamadneh, M. Abu Hammad, and M. Al-Sabbagh, Tubular surfaces of finite type Gauss map, J. Geom. Graph. \textbf{25} (2021), 45-52.
\bibitem{A17} H. Al-Zoubi, F. Abdel-Fattah, M. Al-Sabbagh, Surfaces of finite $III$-type in the Euclidean 3-space, WSEAS Trans. Math. \textbf{20} (2021), 729-735.
\bibitem{B2} Ch. Baikoussis, D. E. Blair, On the Gauss map of Ruled Surfaces, Glasgow Math. J. \textbf{34} (1992), 355-359.
\bibitem{B3} Ch. Baikoussis, F. Denever, P. Emprechts, L. Verstraelen, On the Gauss map of the cyclides of Dupin, Soochow J. Math. \textbf{19} (1993), 417-428.
\bibitem{B4} Ch. Baikoussis, B.-Y. Chen, L. Verstraelen, Ruled Surfaces and tubes with finite type Gauss map, Tokyo J. Math. \textbf{16} (1993), 341-349.
\bibitem{B5} Ch. Baikoussis, L. Verstraelen, On the Gauss map of translation surfaces, Rend. Semi. Mat. Messina Ser II (in press).
\bibitem{B6} Ch. Baikoussis, L. Verstraelen, On the Gauss map of helicoidal surfaces, Rend. Semi. Mat. Messina Ser II \textbf{16} (1993), 31-42.

\bibitem{C3} B.-Y. Chen, Total mean curvature and submanifolds of finite type. Second edition, World Scientific Publisher, (2015).
\bibitem{C5} B.-Y. Chen, Some open problems and conjectures on submanifolds of finite type, Soochow J. Math. \textbf{17} (1991), 169-188.

%
\bibitem{C7} B.-Y. Chen, F. Dillen, L. Verstraelen, L. Vrancken, Ruled surfaces of finite type, Bull. Austral. Math. Soc. \textbf{42} (1990), 447-453.
%
\bibitem{C8} B.-Y. Chen, F. Dillen, Quadrics of finite type, J. of Geom. \textbf{38} (1990), 16-22.
%

\bibitem{D1} F. Defever, R. Deszcz and L. Verstraelen, The compact cyclides of Dupin and a conjecture of B.-Y. Chen. J. Geom. \textbf{46} (1993), 33-38.
%
%
\bibitem{S00} B. Senoussi, H. Al-Zoubi, Translation surfaces of finite type in Sol$_{3}$, Comment. Math. Univ. Carolin. \textbf{61} (2020), 237–256
\bibitem{S0} B. Senoussi and M. Bekkar, Helicoidal surfaces with $\triangle^{J}r = Ar$ in 3-dimensional Euclidean space,
Stud. Univ. Babes-Bolyai Math. \textbf{8} (2015), 437-448.
\bibitem{S1} S. Stamatakis, H. Al-Zoubi, On surfaces of finite Chen-type, Results. Math. \textbf{43} (2003), 181-190.
\bibitem{S2} S. Stamatakis, H. Al-Zoubi, Surfaces of revolution satisfying $\triangle^{III}\mathbf{x}=A\mathbf{x}$, J. Geom. Graph.
\textbf{14} (2010), 181-186.

\end{thebibliography}
\end{document}